\documentclass{amsart}
\usepackage{amssymb}
\usepackage{amscd}

\title
{Pentagon and hexagon equations}
\author{Hidekazu Furusho}

\address{Graduate School of Mathematics, Nagoya University, 
Chikusa-ku, Furo-cho, Nagoya, 464-8602, Japan }
\email{furusho@math.nagoya-u.ac.jp}

\address{D\'{e}partement de math\'{e}matiques et applications,
\'{E}cole Normale Sup\'{e}rieure,
45, rue d'Ulm, F 75230 Paris cedex 05, France}
\email{Hidekazu.Furusho@ens.fr}

\newtheorem{thm}{Theorem}
\newtheorem{lem}[thm]{Lemma}
\newtheorem{cor}[thm]{Corollary}

\theoremstyle{remark}
\newtheorem{ack}{Acknowledgments}        

\theoremstyle{definition}

\newtheorem{rem}[thm]{Remark}

\newtheorem{pf}{Proof}                 
    
\newtheorem{q}[thm]{Question}

\begin{document}
\bibliographystyle{amsalpha+}
\maketitle
\begin{abstract}
The author will prove that Drinfel'd's
pentagon equation implies his two hexagon equations
in the Lie algebra, pro-unipotent, pro-$l$ 
and pro-nilpotent contexts.
\end{abstract}

\tableofcontents

\setcounter{section}{-1}
\section{Introduction}
In his celebrated papers on quantum groups \cite{Dr86, Dr90, Dr} 
Drinfel'd came to the notion of 
quasitriangular quasi-Hopf quantized universal enveloping algebra.
It is a topological algebra which
differs from a topological Hopf algebra 
in the sense that the coassociativity axiom and the cocommutativity
axiom is twisted by an associator and an R-matrix
satisfying a pentagon axiom and two hexagon axioms.
One of the main theorems in \cite{Dr} is that
any quasitriangular quasi-Hopf quantized universal enveloping algebra
modulo twists (in other words gauge transformations \cite{K})
is obtained as a quantization of a pair (called its classical limit) of 
a Lie algebra and its symmetric invariant 2-tensor. 
Quantizations are constructed by universal associators.
The set of group-like universal associators forms a pro-algebraic variety,
denoted $M$.
Its non-emptiness is another of his main theorem
(reproved in \cite{B}).
Our first theorem is on the defining equations of $M$. 

Let us fix notation and conventions:
Let $k$ be a field of characteristic $0$,
$\bar k$ its algebraic closure and
$U\frak F_2=k\langle\langle X,Y\rangle\rangle$ a non-commutative
formal power series ring with two variables $X$ and $Y$.
Its element $\varphi=\varphi(X,Y)$ is called {\it group-like} if it satisfies
$\Delta (\varphi)=\varphi\otimes \varphi$ with $\Delta(X)=X\otimes 1+1\otimes X$ and
$\Delta(Y)=Y\otimes 1+1\otimes Y$ and 
its constant term is equal to $1$.
Its coefficient of $XY$ is denoted by $c_2(\varphi)$.
For any $k$-algebra homomorphism $\iota:U\frak F_2\to S$
the image $\iota(\varphi)\in S$ is denoted 
by $\varphi(\iota(X),\iota(Y))$.
Let $\frak a_4$ be the completion (with respect to the natural grading)
of the Lie algebra over $k$ with generators
$t_{ij}$ ($1\leqslant i,j \leqslant 4$) 
and defining relations $t_{ii}=0$,
$t_{ij}=t_{ji}$, $[t_{ij},t_{ik}+t_{jk}]=0$ ($i$,$j$,$k$: all distinct)
and $[t_{ij},t_{kl}]=0$ ($i$,$j$,$k$,$l$: all distinct).

\begin{thm}\label{M}
Let $\varphi=\varphi(X,Y)$ be a 
group-like element of $U\frak F_2$. 
Suppose that $\varphi$ satisfies Drinfel'd's  pentagon equation:
\begin{equation}\label{pentagon}
\varphi(t_{12},t_{23}+t_{24})
\varphi(t_{13}+t_{23},t_{34})=
\varphi(t_{23},t_{34})
\varphi(t_{12}+t_{13},t_{24}+t_{34})
\varphi(t_{12},t_{23}).
\end{equation}
Then there exists an element (unique up to signature) 
$\mu\in \bar k$
such that the pair $(\mu, \varphi)$ satisfies
his two hexagon equations:
\begin{equation}\label{hexagon}
\exp\{\frac{\mu (t_{13}+t_{23})}{2}\}=
\varphi(t_{13},t_{12})\exp\{\frac{\mu t_{13}}{2}\}
\varphi(t_{13},t_{23})^{-1}
\exp\{\frac{\mu t_{23}}{2}\} \varphi(t_{12},t_{23}),
\end{equation}
\begin{equation}\label{hexagon-b}
\exp\{\frac{\mu (t_{12}+t_{13})}{2}\}=
\varphi(t_{23},t_{13})^{-1}\exp\{\frac{\mu t_{13}}{2}\}
\varphi(t_{12},t_{13})
\exp\{\frac{\mu t_{12}}{2}\} \varphi(t_{12},t_{23})^{-1}.
\end{equation}
Actually this $\mu$ is equal to $\pm(24c_2(\varphi))^\frac{1}{2}$.
\end{thm}

It should be noted that 
we need to use an (actually quadratic) extension  of a field $k$
in order to obtain the hexagon equations from the pentagon equation.
The associator set ${\underline M}$  is 
the pro-algebraic variety 
whose set of $k$-valued points consists of pairs $(\mu,\varphi)$
satisfying \eqref{pentagon}, \eqref{hexagon} and
\eqref{hexagon-b} and 
$M$ is its open subvariety defined by $\mu\neq 0$.
The theorem says that the pentagon equation 
is essentially a single defining equation of the associator set.
The Drinfel'd associator $\varPhi_{\mathrm KZ}\in\bold R
\langle\langle X,Y\rangle\rangle$ 
is a group-like series constructed by solutions of KZ-equation \cite{Dr}.
It satisfies
\eqref{pentagon}, \eqref{hexagon} and \eqref{hexagon-b} with 
$\mu=\pm 2\pi \sqrt{-1}$.
Its coefficients are expressed by multiple zeta values \cite{LM}
(and also \cite{F03}).
The theorem also says that the two hexagon equations do not provide
any new relations 
under the pentagon equation.  

The category of representations of a 
quasitriangular quasi-Hopf quantized universal enveloping algebra
forms a quasitensored category \cite{Dr},
in other words, a braided tensor category \cite{JS};
its associativity constraint and its commutativity constraint
are subject to one pentagon axiom and two hexagon axioms. 
The Grothendieck-Teichm\"{u}ller pro-algebraic group $GT$
is introduced in \cite{Dr} as
a group of deformations of the category which change
its associativity constraint and its commutativity constraint
keeping all three axioms.
It is also conjecturally related to the motivic Galois group of 
${\bold Z}$ (explained in \cite{A}).
Relating to the absolute Galois group $Gal(\bar{\bold Q}/\bold Q)$
of $\bold Q$
its profinite group version $\widehat{GT}$ is discussed in \cite{Ih,S}.
Our second theorem is on defining equations of $GT$. 

\begin{thm}\label{GT}
Let $F_2(k)$ be the the free pro-unipotent algebraic group 
with two variables $x$ and $y$.
Suppose that its element $f=f(x,y)$ satisfies Drinfel'd's pentagon equation:
\begin{equation}\label{GT5}
f(x_{12},x_{23}x_{24})
f(x_{13}x_{23},x_{34})=
f(x_{23},x_{34})
f(x_{12}x_{13},x_{24}x_{34})
f(x_{12},x_{23})
\end{equation}
in $K_4(k)$.
Then there exists an element (unique up to signature)
$\lambda\in \bar k$ such that
the pair $(\lambda,f)$ satisfies 
his hexagon equations (3- and 2-cycle relation):
\begin{equation}\label{GT3}
f(z,x)z^mf(y,z)y^mf(x,y)x^m=1
\text{    with  } xyz=1 \text{ and } m=\frac{\lambda-1}{2},
\end{equation}
\begin{equation}\label{GT2}
f(x,y)f(y,x)=1.
\end{equation}
Actually this $\lambda$ is equal to $\pm(24c_2(f)+1)^\frac{1}{2}$
where $c_2(f)$ stands for $c_2(f(e^X,e^Y))$.
\end{thm}

Here $K_4(k)$ stands for the unipotent completion of 
the pure braid group $K_4=\ker\{B_4\to\frak S_4\}$
of 4 strings
($B_4$: the Artin braid group and 
$\frak S_4$: the symmetric group)
with standard generators $x_{ij}$
($1\leqslant i,j\leqslant 4$).

It should be noted again that 
we need to use an (actually quadratic) extension  of a field $k$
in order to obtain the hexagon equations from the pentagon equation.
The set of pairs $(\lambda,f)$
satisfying \eqref{GT5}, \eqref{GT3} and \eqref{GT2} 
determines a pro-algebraic variety $\underline{GT}$ and
$GT$ is 
its open subvariety defined by $\lambda\neq 0$.
The product structure on $GT(k)$ is given by 
$(\lambda_1,f_1)\circ(\lambda_2,f_2):=(\lambda,f)$
with $\lambda=\lambda_1\lambda_2$ and
$f(x,y)=f_1(f_2x^{\lambda_2}f^{-1}_2,y^{\lambda_2})f_2$.
The theorem says that the pentagon equation 
is essentially a single defining equation of $GT$.
 
The construction of the paper is as follows:
\S\ref{Lie algebra case} is a crucial part of the paper.
The implication of the pentagon equation is proved for Lie series.
In \S\ref{proof of main theorem} 
we give a proof of theorem \ref{M} by using Drinfel'd's gadgets.
\S\ref{proof of theorem GT} gives a proof of
theorem \ref{GT} and its analogue in the pro-$l$ group and
pro-nilpotent group setting.

\section{Lie algebra case}\label{Lie algebra case}
In this section we prove the Lie algebra version of Theorem \ref{M}
in a rather combinatorial argument.

Let $\frak F_2$ be the set of 
Lie-like elements $\varphi$ in $U\frak F_2$
(i.e. $\Delta(\varphi)=\varphi\otimes 1+1\otimes \varphi$).

\begin{thm}\label{Lie}
Let $\varphi$ be a commutator Lie-like element
\footnote{
In this paper we call a series $\varphi=\varphi(X,Y)$ {\it commutator Lie-like}
if it is Lie-like and
its coefficient of $X$ and $Y$ are both $0$,
in other words $\varphi\in \frak F'_2:=[\frak F_2,\frak F_2]$.
}
with $c_2(\varphi)=0$.
Suppose that $\varphi$ satisfies the pentagon equation (5-cycle relation):
\begin{equation}\label{SDA5}
 \varphi(X_{12},X_{23})
+\varphi(X_{34},X_{45})
+\varphi(X_{51},X_{12})
+\varphi(X_{23},X_{34})
+\varphi(X_{45},X_{51})=0
\end{equation} 
in $\hat{\frak P_5}$.
Then it also satisfies the hexagon equations 
(3- and 2-cycle relation):
\begin{equation}\label{SDA3}
\varphi(X,Y)+\varphi(Y,Z)+\varphi(Z,X)=0 \text{ with } X+Y+Z=0,
\end{equation}
\begin{equation}\label{SDA2}
\varphi(X,Y)+\varphi(Y,X)=0.
\end{equation}
\end{thm}

Here $\hat{\frak P_5}$ stands for the completion 
(with respect to the natural grading) of 
the pure sphere braid Lie algebra
$\frak P_5$ \cite{Ih} with 5 strings;
the Lie algebra 
generated by $X_{ij}$ ($1\leqslant i,j\leqslant 5$) with clear relations
$X_{ii}=0$, $X_{ij}=X_{ji}$, $\sum_{j=1}^5 X_{ij}=0$
($1\leqslant i,j\leqslant 5$) and
$[X_{ij},X_{kl}]=0$ if $\{i,j\}\cap\{k,l\}=\emptyset$.
It is a quotient of $\frak a_4$ (cf. \S\ref{proof of main theorem}).

\begin{pf}
There is a projection from $\hat{\frak P_5}$ to  
the completed free Lie algebra $\frak F_2$ generated by $X$ and $Y$ 
by putting $X_{i5}=0$, $X_{12}=X$ and $X_{23}=Y$.
The image of $5$-cycle relation gives $2$-cycle relation.

For our convenience we denote 
$\varphi(X_{ij},X_{jk})$ 
($1\leqslant i,j,k \leqslant 5$)
by $\varphi_{ijk}$.
Then the $5$-cycle relation can be read as
$$
\varphi_{123}+\varphi_{345}+\varphi_{512}+\varphi_{234}+\varphi_{451}=0.
$$
We denote LHS by $P$.
Put $\sigma_i$ ($1\leqslant i \leqslant 4$)
to be elements of $\frak S_5$ defined as follows:
$\sigma_1(12345)=(12345)$,
$\sigma_2(12345)=(54231)$,
$\sigma_3(12345)=(13425)$ and
$\sigma_4(12345)=(43125)$.
Then 
\begin{align*}
\sum\nolimits_{i=1}^{4}\sigma_i(P)=
&\quad\varphi_{123}+\varphi_{345}+\varphi_{512}+\varphi_{234}+\varphi_{451} \\
&+\varphi_{542}+\varphi_{231}+\varphi_{154}+\varphi_{423}+\varphi_{315}\\
&+\varphi_{134}+\varphi_{425}+\varphi_{513}+\varphi_{342}+\varphi_{251}\\
&+\varphi_{431}+\varphi_{125}+\varphi_{543}+\varphi_{312}+\varphi_{254}.\\
\end{align*}
By the 2-cycle relation,
$\varphi_{ijk}=-\varphi_{kji}$ ($1\leqslant i,j,k \leqslant 5$).
This gives
\begin{align*}
\sum\nolimits_{i=1}^{4}\sigma_i(P)=&
\quad(\varphi_{123}+\varphi_{231}+\varphi_{312})
+(\varphi_{512}+\varphi_{125}+\varphi_{251})\\
&+(\varphi_{234}+\varphi_{342}+\varphi_{423})
+(\varphi_{542}+\varphi_{425}+\varphi_{254}).
\end{align*}

By $[X_{23},X_{12}+X_{23}+X_{31}]
=[X_{31},X_{12}+X_{23}+X_{31}]
=[X_{12},X_{12}+X_{23}+X_{31}]=0$ 
and $\varphi\in\frak F'_2$,
$\varphi_{231}=\varphi(X_{23},X_{31})=\varphi(X_{23},-X_{12}-X_{23})$ and
$\varphi_{312}=\varphi(X_{31},X_{12})=\varphi(-X_{12}-X_{23},X_{12})$.

By $[X_{51},X_{12}+X_{25}+X_{51}]
=[X_{12},X_{12}+X_{25}+X_{51}]
=[X_{25},X_{12}+X_{25}+X_{51}]=0$ 
and $\varphi\in\frak F'_2$,
$\varphi_{512}=\varphi(X_{51},X_{12})=\varphi(-X_{12}-X_{25},X_{12})$ and
$\varphi_{251}=\varphi(X_{25},X_{51})=\varphi(X_{25},-X_{12}-X_{25})$.

By $[X_{23},X_{42}+X_{23}+X_{34}]
=[X_{34},X_{42}+X_{23}+X_{34}]
=[X_{42},X_{42}+X_{23}+X_{34}]=0$ 
and $\varphi\in\frak F'_2$,
$\varphi_{234}=\varphi(X_{23},X_{34})=\varphi(X_{23},-X_{42}-X_{23})$ and
$\varphi_{342}=\varphi(X_{34},X_{42})=\varphi(-X_{42}-X_{23},X_{42})$.

By $[X_{54},X_{42}+X_{25}+X_{54}]
=[X_{42},X_{42}+X_{25}+X_{54}]
=[X_{25},X_{42}+X_{25}+X_{54}]=0$ 
and $\varphi\in\frak F'_2$,
$\varphi_{542}=\varphi(X_{54},X_{42})=\varphi(-X_{42}-X_{25},X_{42})$ and
$\varphi_{254}=\varphi(X_{25},X_{54})=\varphi(X_{25},-X_{42}-X_{25})$.

Put $R(X,Y)=\varphi(X,Y)+\varphi(Y,-X-Y)+\varphi(-X-Y,X)$.
Then 
$$
\sum\nolimits_{i=1}^{4}\sigma_i(P)=
R(X_{21},X_{23})
+R(X_{21},X_{25})
+R(X_{24},X_{23})
+R(X_{24},X_{25}).
$$

The elements $X_{21}$, $X_{23}$, $X_{24}$ and $X_{25}$ 
generate a completed Lie subalgebra $\frak F_3$ of $\hat{\frak P_5}$ 
which is free of rank $3$ and whose  set of relations is given by 
$X_{21}+X_{23}+X_{24}+X_{25}=0$.
It contains $\sum\nolimits_{i=1}^{4}\sigma_i(P)$.
Let $q_1:\frak F_3\to\frak F_2$ be the projection sending
$X_{21}\mapsto X$, $X_{23}\mapsto Y$ and $X_{24}\mapsto X$.
Then
$$
q_1(\sum\nolimits_{i=1}^{4}\sigma_i(P))=
R(X,Y)+R(X,-2X-Y)+R(X,Y)+R(X,-2X-Y).
$$
Since $P=0$, we have 
$R(X,-2X-Y)=-R(X,Y)$.
Let $q_2:\frak F_3\to\frak F_2$ be the projection sending
$X_{21}\mapsto X$, $X_{23}\mapsto X$ and $X_{24}\mapsto Y$.
Then
$$
q_2(\sum\nolimits_{i=1}^{4}\sigma_i(P))=
R(X,X)+R(X,-2X-Y)+R(Y,X)+R(Y,-2X-Y).
$$
By $\varphi\in\frak F'_2$, $R(X,X)=0$.
By definition,
$R(Y,-2X-Y)=R(2X,Y)$.
Since $P=0$, 
$-R(X,Y)+R(Y,X)+R(Y,2X)=0$.
The $2$-cycle relation gives $R(X,Y)=-R(Y,X)$. 
Therefore
$2R(X,Y)=R(2X,Y)$.
Expanding this equation in terms of a linear basis, say the Hall basis,
we see that $R(X,Y)$ must be of the form
$\sum\nolimits_{m=1}^\infty a_m(adY)^{m-1}(X)$ with $a_m\in k$.
Since it satisfies $R(X,Y)=-R(Y,X)$,
we have $a_1=a_3=a_4=a_5=\cdots=0$.
By our assumption $c_2(\varphi)=0$,
$a_2$ must be $0$ also.
Therefore $R(X,Y)=0$,
which is the $3$-cycle relation.
\qed
\end{pf}

We note that the assumption $c_2(\varphi)=0$ is necessary:
e.g. the element $\varphi=[X,Y]$ satisfies the $5$-cycle relation but 
it does not satisfy the $3$-cycle relation.

\begin{rem}
There is partially a geometric picture in the proof:
We have a de Rham fundamental groupoid \cite{De} (see also \cite{F07}) 
of the moduli
${\mathcal M}_{0,n}=
\{(x_1:\cdots:x_n)\in (\bold P^1)^n|x_i\neq x_j (j\neq j)\}/PGL(2)$
for $n\geqslant 4$,
its central extension given by the normal bundle of $\mathcal M_{0,n-1}$
inside its stable compactification $\overline{\mathcal M}_{0,n}$
and maps between them.
An automorphism of the system is determined by considering 
what happens to the canonical de Rham path from `$0$' to `$1$' (loc.cit.)
in $\mathcal M_{0,4}=\bold P^1\backslash\{0,1,\infty\}$.
Equation \eqref{SDA5} reflects the necessary condition
that such an automorphism must keep the property
that the image of the composite of the path,  
the boundaries of the fundamental pentagon $\mathcal B_5$ \cite{Ih}
formed by the divisors $x_i=x_{i+1}$ ($i\in\bold Z/5\bold Z$) in 
$\overline{\mathcal M}_{0,5}(\bold R)$,
must be a trivial loop.
Each $\sigma_i(\mathcal B_5)$ ($1\leqslant i \leqslant 4$) 
is a connected component of $\mathcal M_{0,5}(\bold R)$. 
The sum of four 5-cycles
$\sum_{i=1}^{4}\sigma_i(P)$ corresponds to a path following 
the (oriented) boundaries of the four real pentagonal regions
$\sigma_i(\mathcal B_5)$ of $\mathcal M_{0,5}(\bold R)$. 
The four 3-cycles  
correspond to four loops around the four boundary divisors 
$x_{4}=x_{5}$, $x_{3}=x_{4}$, $x_{5}=x_{1}$ and $x_{1}=x_{3} $ in 
$\overline{\mathcal M}_{0,5}(\bold R)$. 
The author expects that the geometric interpretation might help
to adapt our proof to the pro-finite context
(cf. question \ref{pro-finite}).
\end{rem}

The equations  \eqref{SDA5},\eqref{SDA3} and \eqref{SDA2}
are defining equations of Ihara's  
stable derivation (Lie-)algebra \cite{Ih}.
Its Lie bracket is given by 
$<\varphi_1,\varphi_2>:=[\varphi_1,\varphi_2]
+D_{\varphi_2}(\varphi_1)-D_{\varphi_1}(\varphi_2)$
where $D_\varphi$ is the derivation of $\frak F_2$
given by $D_\varphi(X)=[\varphi,X]$ and $D_\varphi(B)=0$.
We note that its completion with respect to degree is equal to 
the graded Lie algebra $\frak{grt}_1$ 
of the Grothendieck-Teichm\"{u}ller group $GT$ in \cite{Dr}.
Our theorem says that
the pentagon equation is its single defining equation 
and two hexagon equations are needless for its definition
when $\deg \varphi\geqslant 3$.

\section{Proof of Theorem \ref{M}}\label{proof of main theorem}
This section is devoted to a proof of Theorem \ref{M}.

Between the Lie algebra $\frak a_4$ in theorem \ref{M} 
and $\hat{\frak P_5}$ in theorem \ref{Lie}
there is a natural surjection $\tau:\frak a_4\to\hat{\frak P_5}$
sending $t_{ij}$ to $X_{ij}$ ($1\leqslant i,j\leqslant 4$).
Its kernel is generated by 
$\Omega=\sum_{1\leqslant i<j\leqslant 4}t_{ij}$.
We also denote its induced morphism 
$U\frak a_4\to U\hat{\frak P_5}$ by $\tau$.
On the pentagon equation we have

\begin{lem}\label{pentagon=5}
Let $\varphi$ be a group-like element.
Giving the pentagon equation \eqref{pentagon} for $\varphi$ 
is equivalent to giving that
$\varphi$ is commutator group-like 
\footnote{
In this paper we call a series $\varphi=\varphi(X,Y)$ {\it commutator group-like}
if it is group-like and
its coefficient of $X$ and $Y$ are both $0$.}
and $\varphi$ satisfies the 5-cycle relation in $U\hat{\frak P_5}$: 
\begin{equation}\label{5-cycle}
 \varphi(X_{12},X_{23})
\varphi(X_{34},X_{45})
\varphi(X_{51},X_{12})
\varphi(X_{23},X_{34})
\varphi(X_{45},X_{51})=1.
\end{equation}
\end{lem}

\begin{pf}
Assume \eqref{pentagon}. 
Denote the abelianization of 
$\varphi(X,Y)\in k\langle\langle X,Y\rangle\rangle$
by $\varphi^{\mathrm{ab}}\in k[[X,Y]]$.
The series  $\varphi$ is group-like, 
so $\varphi^{\mathrm{ab}}$ is also,
i.e. $\Delta(\varphi^{\mathrm{ab}})=
\varphi^{\mathrm{ab}}\otimes\varphi^{\mathrm{ab}}$.
Therefore $\varphi^{\mathrm{ab}}$
must be of the form $\exp\{\alpha X+\beta Y\}$ with $\alpha,\beta\in k$.
The equation  \eqref{pentagon} gives $\alpha X_{12}+\beta X_{34}=0$.
Hence $\alpha=\beta=0$ which means that $\varphi$ is commutator group-like.
Therefore \par
$\varphi(X_{12},X_{51})=\varphi(X_{12},-X_{21}-X_{52})
=\varphi(X_{12},X_{23}+X_{24})$\\
by $[X_{12},X_{51}+X_{21}+X_{52}]=[X_{51},X_{51}+X_{21}+X_{52}]=0$,\par
$\varphi(X_{45},X_{34})=\varphi(-X_{43}-X_{53},X_{34})
=\varphi(X_{13}+X_{23},X_{34})$\\
by $[X_{45},X_{45}+X_{43}+X_{53}]=[X_{34},X_{45}+X_{43}+X_{53}]=0$
and \par
$\varphi(X_{45},X_{51})=\varphi(-X_{14}-X_{15},X_{51})
=\varphi(-X_{14}-X_{15},-X_{14}-X_{45})
=\varphi(X_{12}+X_{13},X_{24}+X_{34})$\\
by $[X_{45},X_{45}+X_{14}+X_{51}]=[X_{51},X_{45}+X_{14}+X_{51}]=0$
and $[X_{14}+X_{15},X_{51}+X_{14}+X_{45}]
=[X_{51},X_{51}+X_{14}+X_{45}]=0$.
(N.B. If $\varphi$ is commutator group-like,
$\varphi(A+C,B)=\varphi(A,B+C)=\varphi(A,B)$
with $[A,C]=[B,C]=0$.)
So the image of \eqref{pentagon} by $\tau$ is
\begin{equation}\label{gokakukei}
\varphi(X_{12},X_{51})\varphi(X_{45},X_{34})=
\varphi(X_{23},X_{34})\varphi(X_{45},X_{51})\varphi(X_{12},X_{23}).
\end{equation}
The following lemma below gives \eqref{5-cycle}.

Conversely assume \eqref{5-cycle} and the commutator group-likeness for $\varphi$.
The lemma below gives the above equality \eqref{gokakukei}.
Whence we say \eqref{pentagon} modulo $\ker\tau$.
That is, the quotient of LHS of \eqref{pentagon} by RHS of \eqref{pentagon}
is expressed as $\exp \gamma\Omega$ for some $\gamma\in k$.
Since both hand sides of \eqref{pentagon} are commutator group-like,  
$\exp \gamma\Omega$ must be also.
Therefore $\gamma$ must be $0$, which gives \eqref{pentagon}.
\qed
\end{pf}

\begin{lem}
Let $\varphi$ be a group-like element. 
If $\varphi$ is commutator group-like and 
it satisfies the $5$-cycle relation \eqref{5-cycle},
it also satisfies the $2$-cycle relation:
\begin{equation}\label{2-cycle}
\varphi(X,Y)\varphi(Y,X)=1.
\end{equation}

And if $\varphi$ satisfies the pentagon equation \eqref{pentagon},
it also satisfies \eqref{2-cycle}.
\end{lem}

\begin{pf}
There is a projection $U\frak P_5^\land\to U\frak F_2$ 
by putting $X_{i5}=0$ ($1\leqslant i\leqslant 5$), 
$X_{12}=X$ and $X_{23}=Y$.
The image of \eqref{5-cycle} is \eqref{2-cycle}
by the commutator group-likeness.

As was shown in the previous lemma,
the equation \eqref{pentagon} for $\varphi$ in $U\frak a_4$ 
implies its commutator group-likeness and
\eqref{gokakukei} in $U\hat{\frak P_5}$.
The image of \eqref{gokakukei} by the projection
gives \eqref{2-cycle}.
\qed
\end{pf}

In \cite{IM} they showed the equivalence 
between \eqref{pentagon} and \eqref{5-cycle}
assuming the commutatativity and the 2-cycle relation
in the pro-finite group setting.
But by the above argument the latter assumption can be excluded.

As for the hexagon equations we also have

\begin{lem}
Let $\varphi$ be a group-like element.
Giving two hexagon equations \eqref{hexagon} and \eqref{hexagon-b} 
for $\varphi$ 
is equivalent to giving the 2-cycle relation \eqref{2-cycle} and 
the 3-cycle relation:
\begin{equation}\label{3-cycle}
e^{\frac{\mu X}{2}}\varphi(Z,X)e^{\frac{\mu Z}{2}}\varphi(Y,Z)
e^{\frac{\mu Y}{2}} \varphi(X,Y)=1
\text{    with  } X+Y+Z=0.
\end{equation}
\end{lem}

\begin{pf}
We review the proof in \cite{Dr}.
The Lie subalgebra generated by $t_{12}$, $t_{13}$ and $t_{23}$
is  the direct sum of its center, generated by $t_{12}+t_{23}+t_{13}$,
and the free Lie algebra generated by $X=t_{12}$ and $Y=t_{23}$.
The projections of \eqref{hexagon} and \eqref{hexagon-b} 
to the first component are both tautologies but
the projections to the second component are
$
e^{\frac{\mu X}{2}}\varphi(Z,X)e^{\frac{\mu Z}{2}}\varphi(Z,Y)^{-1}
e^{\frac{\mu Y}{2}} \varphi(X,Y)=1
$ and
$
e^{\frac{\mu X}{2}}\varphi(Z,X)e^{\frac{\mu Z}{2}}\varphi(Z,Y)^{-1}
e^{\frac{\mu Y}{2}} \varphi(Y,X)^{-1}=1$.
They are equivalent to \eqref{2-cycle} and \eqref{3-cycle}.
\qed
\end{pf}

The followings are keys to prove theorem \ref{M}.

\begin{lem}\label{special}
Let $\varphi_1$ and $\varphi_2$ be
commutator group-like elements.
Put $\varphi_3=\varphi_2\circ\varphi_1(X,Y):=
\varphi_2(\varphi_1X\varphi_1^{-1},Y)\cdot \varphi_1$.
Assume that $\varphi_1$ satisfies 
\eqref{5-cycle}, \eqref{2-cycle} and 
\begin{equation}\label{special 3-cycle}
\varphi(Z,X)\varphi(Y,Z)\varphi(X,Y)=1
\text{    with  } X+Y+Z=0.
\end{equation}
Then $\varphi_2$ satisfies \eqref{5-cycle}
if and only if 
$\varphi_3$ satisfies \eqref{5-cycle}.
\end{lem}

\begin{pf}
By the arguments in \cite{S}\S1.2,
$\varphi_1$ determines an automorphism of $U\frak P_5^\land$
sending 
$X_{12}\mapsto X_{12}$,
$X_{23}\mapsto \varphi_1(X_{12},X_{23})^{-1} X_{23}\varphi_1(X_{12},X_{23})$,
$X_{34}\mapsto \varphi_1(X_{34},X_{45}) \\
X_{34} \varphi_1(X_{34},X_{45})^{-1}$,
$X_{45}\mapsto X_{45}$ and 
$X_{51}\mapsto \varphi_1(X_{12},X_{23})^{-1} 
\varphi_1(X_{45},X_{51})^{-1}
X_{51} \\
\varphi_1(X_{45},X_{51})
\varphi_1(X_{12},X_{23})$.
The direct calculation shows that LHS of \eqref{5-cycle} for $\varphi_2$
maps to LHS of \eqref{5-cycle} for $\varphi_3(X,Y)$.
This gives the claim. 
\qed
\end{pf}

\begin{lem}
Let $\varphi$ be a commutator group-like element with $c_2(\varphi)=0$. 
Suppose that $\varphi$ satisfies \eqref{5-cycle}.
Then it also satisfies
\eqref{special 3-cycle}.
\end{lem}

\begin{pf}
The proof is given by induction.
Suppose that we have \eqref{special 3-cycle}$\mod \deg n$.
The element $\varphi$ satisfies the commutator group-likeness,
\eqref{5-cycle}, \eqref{2-cycle} and 
\eqref{special 3-cycle}$\mod \deg n$, in other words,
it is an element of algebraic group $GRT_1^{(n)}(k)$ \cite{Dr}\S 5.
Denote its corresponding Lie element by $\psi$.
It is an element of the Lie algebra $\frak{grt}_1^{(n)}(k)$ 
(loc.cit), that means,
it is expressed by
$\psi=\sum\nolimits_{i=3}^{n-1}\psi^{(i)}\in
k\langle\langle X,Y\rangle\rangle$ 
where $\psi^{(i)}$ is a homogeneous Lie element
with $\deg \psi^{(i)}=i$ and
satisfies \eqref{SDA5}, \eqref{SDA3} and \eqref{SDA2}.
The Lie algebra $\frak{grt}_1(k)=\varprojlim\frak{grt}_1^{(n)}(k)$
is graded by degree and $\psi$ also determines an element 
(denoted by the same symbol $\psi$)
of $\frak{grt}_1(k)$.
Let $\text{Exp}:\frak{grt}_1(k)
\to GRT_1(k)=\varprojlim GRT_1^{(n)}(k)$
be the exponential morphism.
Put $\varphi_1=\text{Exp }\psi$.
It is commutator group-like and
it satisfies \eqref{5-cycle},
\eqref{2-cycle}, \eqref{special 3-cycle} and
$\varphi\equiv\varphi_1\mod\deg n$ (loc.cit).
Let $\varphi_2$ be a series defined by  
$\varphi=\varphi_2\circ\varphi_1$.
Then $\varphi_2$ is commutator group-like and
it satisfies \eqref{5-cycle} by lemma \ref{special}.
By $\varphi\equiv\varphi_1\mod\deg n$,
$\varphi_2\equiv 1\mod\deg n$.
Denote the degree $n$-part of $\varphi_2$ by $\psi^{(n)}$.
Because $\varphi_2\equiv 1+\psi^{(n)}\mod\deg n+1$,
\eqref{5-cycle} for $\varphi_2$ yields \eqref{SDA5} for $\psi^{(n)}$ and
the group-likeness for  $\varphi_2$ yields the Lie-likeness for $\psi^{(n)}$.
By theorem \ref{Lie}, $\psi^{(n)}$ satisfies \eqref{SDA3} and \eqref{SDA2},
which means $\psi^{(n)}\in\frak{grt}_1(k)$.
Since $\text{Exp }\psi^{(n)}\in GRT_1(k)$ and
$\varphi_2\equiv \text{Exp }\psi^{(n)}\mod\deg n+1$,
$\varphi_2$ belongs to $GRT_1^{(n+1)}(k)$.
Since $\varphi_1$ also determines an element of $GRT_1^{(n+1)}(k)$, 
$\varphi$ must belong to $GRT_1^{(n+1)}(k)$.
This means that $\varphi$ satisfies \eqref{special 3-cycle}$\mod \deg n+1$.
\qed
\end{pf}

\begin{thm}
Let $\varphi$ be a commutator group-like element. 
Suppose that $\varphi$ satisfies the $5$-cycle relation \eqref{5-cycle}.
Then there exists an element (unique up to signature) 
$\mu\in \bar k$
such that the pair $(\mu, \varphi)$ satisfies
the $3$-cycle relation \eqref{3-cycle}.
Actually this $\mu$ is equal to $\pm(24c_2(\varphi))^\frac{1}{2}$.
\end{thm}

{\bf Proof.}
We may assume $c_2(\varphi)\neq 0$
by the previous lemma.
Let $\mu$ be a solution of $x^2=24c_2(\varphi)$ in $\bar k^\times$.
Put $M'_\mu$ (resp. $M_\mu$ \cite{Dr}) to
be the pro-affine algebraic variety whose $\bar k$-valued points
are commutator group-like series $\varphi$ in 
$\bar k\langle\langle X,Y\rangle\rangle$ 
satisfying \eqref{5-cycle} and $c_2(\varphi)=\frac{\mu^2}{24}$
(resp. \eqref{5-cycle}, \eqref{2-cycle} and \eqref{3-cycle})
for $(\mu,\varphi)$.
By calculating the coefficient of $XY$
in \eqref{3-cycle} for $(\mu,\varphi)$,
we get $3c_2(\varphi)-\frac{\mu^2}{8}=0$.
Thus $M_\mu$ is a pro-subvariety of $M'_\mu$.
To prove $M'_\mu=M_\mu$,
it suffices to show this for $\mu=1$ because we have
a replacement $\varphi(A,B)$ by $\varphi(\frac{A}{\mu},\frac{B}{\mu})$.
In a similar way to \cite{F06}\S 6
the regular function ring 
$\mathcal{O}(M'_1)$ 
(resp. $\mathcal{O}(M_1)$) is encoded
the weight filtration 
$W=\{W_n\mathcal{O}(M'_1)\}_{n\in\bold Z}$
(resp. $\{W_n\mathcal{O}(M_1)\}_{n\in\bold Z}$):
the algebra $\mathcal{O}(M'_1)$
(resp. $\mathcal{O}(M_1)$) is generated by
$x_W$'s ($W$: word 
\footnote{
A {\it word} means a monic monomial element but $1$
in $k\langle\langle X,Y\rangle\rangle$.
})
and defined by 
the commutator group-likeness,
\eqref{5-cycle} and $c_2(\varphi)=\frac{1}{24}$
(resp. \eqref{5-cycle}, \eqref{2-cycle} and \eqref{3-cycle})
for $\varphi=1+\sum_W x_WW$.
Put $\deg x_W=\deg W$.
Each $W_n\mathcal{O}(M'_1)$
(resp. $W_n\mathcal{O}(M_1)$)
is the vector space generated by 
polynomials whose total degree is less than or equal to $n$.

The inclusion $M_1 \to M'_1$ gives 
a projection  $\mathcal{O}(M'_1)\twoheadrightarrow\mathcal{O}(M_1)$
which is strictly compatible with the filtrations.
It induces a projection
$p:Gr^W_\cdot\mathcal{O}(M'_1)\twoheadrightarrow Gr^W_\cdot\mathcal{O}(M_1)$
between their associated graded quotients.
The graded quotient $Gr^W_\cdot\mathcal{O}(M_1)$
is isomorphic to $\mathcal{O}(GRT_1)$ by \cite{F06} theorem 6.2.2.
It is the algebra generated by $\bar x_W$'s and defined by
the commutator group-likeness,
\eqref{5-cycle}, \eqref{2-cycle} and \eqref{special 3-cycle} for
$\bar\varphi=1+\sum_W \bar x_WW$. 
On the other hand
the graded quotient $Gr^W_\cdot\mathcal{O}(M'_1)$
is generated by $\bar x_W$'s. 
These generators especially satisfy 
the commutator group-likeness, \eqref{5-cycle} 
and $c_2(\bar\varphi)=0$
for  $\bar\varphi=1+\sum_W \bar x_WW$ among others.
By the previous lemmas, $\bar\varphi$ must also satisfy
\eqref{2-cycle} and \eqref{special 3-cycle}.
Therefore $p$ should be an isomorphism.
This implies $M'_1=M_1$.
\qed

The combination of this theorem with the previous lemmas
completes the proof of theorem \ref{M}.

\section{Proof of Theorem \ref{GT}}\label{proof of theorem GT}
In this section we deduce theorem \ref{GT} from the previous theorem
and also show its pro-$l$ group analogue 
(corollary \ref{pro-l}) and its pro-nilpotent group analogue 
(corollary \ref{pro-nilpotent}).

{\bf Proof of theorem \ref{GT}.}
Let $f$ be an element of $F_2(k)$ satisfying \eqref{GT5}.
Let $\lambda$ be a solution of $\frac{x^2-1}{24}=c_2(f)$.
Let $\mu\in k^\times$ and $\varphi\in 
k\langle\langle A,B\rangle\rangle$
be a pair such that $\varphi$ is commutator group-like and
$(\mu,\varphi)$ satisfies
\eqref{5-cycle}, \eqref{2-cycle} and \eqref{3-cycle}.
Put
$\varphi'=f(\varphi e^{\mu X}\varphi^{-1},e^{\mu Y})\cdot \varphi
\in\bar k\langle\langle A,B\rangle\rangle$.
In the proof of \cite{Dr} proposition 5.1
it is shown that 
giving \eqref{GT5} for $f$ is equivalent to 
giving \eqref{pentagon} for $\varphi'$.
Hence $\varphi'$ satisfies \eqref{5-cycle} by lemma \ref{pentagon=5}.
Put $\mu'=\lambda\mu$.
The equation \eqref{3-cycle} for $(\mu,\varphi)$ gives 
$c_2(\varphi)=\frac{\mu^2}{24}$.
So $c_2(\varphi')=c_2(\varphi)+\mu^2c_2(f)=\frac{\mu'^2}{24}$.
Since $\varphi'$ satisfies \eqref{5-cycle},
our previous theorem gives \eqref{3-cycle} for $(\mu',\varphi')$.
Consider the group isomorphism from $F_2(k)$ to the set of 
group-like elements of $U\frak F_2$
which sends $x$ to $e^{\mu X}$ and $y$ to
$e^{-\frac{\mu}{2}X}\varphi(Y,X)e^{\mu Y}\varphi(Y,X)^{-1}
e^{\frac{\mu}{2}X}$.
Consequently $z$ goes to $\varphi(Z,X)e^{\mu Z}\varphi(Z,X)^{-1}$
by \eqref{2-cycle} and \eqref{3-cycle} for $(\mu,\varphi)$. 
The direct calculation shows that 
LHS of \eqref{GT3} maps to LHS of \eqref{3-cycle}.
Therefore
giving \eqref{GT3} for $(\lambda, f)$
is equivalent to
giving \eqref{3-cycle} for $(\mu',\varphi')$. 
This completes the proof of theorem \ref{GT}. 
\qed

\begin{rem}
By the same argument as lemma \ref{pentagon=5},
giving the pentagon equation \eqref{GT5} for $f$ 
is equivalent to giving that
$f(e^X,e^Y)$ is commutator group-like 
and $f$ satisfies the 5-cycle relation in $P_5(k)$: 
$$
f(x_{12},x_{23})
f(x_{34},x_{45})
f(x_{51},x_{12})
f(x_{23},x_{34})
f(x_{45},x_{51})=1.
$$
Here $P_5(k)$ means the unipotent completion 
of the pure sphere braid group
with 5 strings and $x_{ij}$ means its standard generator.
Occasionally, in some of the literature, 
the formula is used directly instead of \eqref{GT5}
in the definition of the Grothendieck-Teichm\"{u}ller group.
\end{rem}

As a corollary the following pro-$l$ ($l$: a prime) group and pro-nilpotent group
version of theorem \ref{GT} 
are obtained by the natural embedding from the pro-$l$ completion 
$F_2^{(l)}$ to $F_2(\bold Q_l)$ and its associated embedding
from the pro-nilpotent completion 
$F_2^{nil}=\prod_{l:\text{a prime}} F_2^{(l)}$ to
$\prod_l F_2(\bold Q_l)$.

\begin{cor}\label{pro-l}
Let $f=f(x,y)$ be an element of $F_2^{(l)}$
satisfying \eqref{GT5} in $K_4^{(l)}$
(: the pro-$l$ completion of $K_4$).
Then there exists $\lambda$ 
such that the pair $(\lambda, f)$ satisfies \eqref{GT3} and \eqref{GT2}.
Actually this $\lambda$ is equal to $\pm(24c_2(f)+1)^\frac{1}{2}$.
\end{cor}

\begin{cor}\label{pro-nilpotent}
Let $f=f(x,y)$ be an element of $F_2^{nil}$
satisfying \eqref{GT5} in $K_4^{nil}=\prod_l K_4^{(l)}$.
Then there exists $\lambda$ 
such that the pair $(\lambda, f)$ satisfies \eqref{GT3} and \eqref{GT2}.
Actually this $\lambda$ is equal to $\pm(24c_2(f)+1)^\frac{1}{2}$.
\end{cor}

It should be noted that though $\lambda$ might lie on a quadratic extension
the equation \eqref{GT3} makes sense for such $(\lambda,f)$.
In the pro-unipotent context taking a quadratic extension is necessary:
the Drinfel'd associator $\varPhi_{\mathrm KZ}\in\bold R
\langle\langle X,Y\rangle\rangle$  satisfies 
\eqref{hexagon} and \eqref{hexagon-b} with 
$\mu=\pm 2\pi \sqrt{-1}\not\in \bold R^\times$.
In the pro-$l$ context the author thinks that
it might also happen
$\pm(24c_2(f)+1)^\frac{1}{2}\not\in \bold Z_l^\times$.

We have a group theoretical definition of $c_2(f)$
(cf. \cite{LS} lemma 9):
put $F^{(l)}_2(1):=[F^{(l)}_2,F^{(l)}_2]$ and 
$F^{(l)}_2(2):=[F^{(l)}_2(1),F^{(l)}_2(1)]$
where $[\cdot,\cdot]$ means the topological commutator.
The quotient group $F^{(l)}_2(1)/ F^{(l)}_2(2)$ is
cyclic generated by the commutator $(x,y)$.
For $f\in F^{(l)}_2(1)$,  $c_2(f)\in \bold Z_l$ is defined by
$(x,y)^{c_2(f)}\equiv f$ in this quotient. 
Posing the following question on 
a pro-finite group analogue
of theorem \ref{GT}
might be particularly interesting:

\begin{q}\label{pro-finite}
Let $f=f(x,y)$ be an element of the pro-finite completion $\hat{F_2}$
satisfying \eqref{GT5} (hence \eqref{GT2}) in the pro-finite completion $\hat{K_4}$.
Let $c_2(f)$ be an element in $\widehat{\bold Z}$ 
defined in a similar way to the above.
Assume that there exists $\lambda$ in $\widehat{\bold Z}$ 
such that $\lambda^2=24c_2(f)+1$.
Then does 
the pair $(\lambda, f)$ satisfy \eqref{GT3}?
\end{q}

\begin{rem}
Although the pentagon equation \eqref{GT5} 
implies the two hexagon equations \eqref{GT3} and \eqref{GT2}
of $GT$,
it does not mean that the pentagon axiom \cite{Dr}(1.7)
implies two hexagon axioms, \cite{Dr}(1.9a) and (1.9b),
of braided tensor categories.
It is because that the pentagon equation \eqref{GT5} of $GT$
is a consequence of the three axioms 
of braided tensor categories:
$GT$ is interpreted as a group of deformations of braided tensor categories
by Drinfel'd in \cite{Dr}\S 4.
The equation \eqref{GT5} of $GT$ is read as a condition to
keep the pentagon axiom. 
However it is formulated in terms of the braid group $K_4$, 
where its generators $x_{ij}$'s are subject to the braid relations.
In his interpretation
the relations are guaranteed by
the dodecagon diagram (the Yang-Baxter equation) 
(see \cite{JS} proposition 2.1 and \cite{K} theorem XIII.1.3)
which is deduced from 
two hexagon axioms. 
\end{rem}

\begin{ack}
The author is grateful to Pierre Deligne
for his crucial comments on the earlier version of the paper.
He particularly thanks Leila Schneps for shortening our proof of theorem \ref{Lie},
Pierre Lochack and referee for valuable comments.
He is supported by Research Aid of Inoue Foundation for Science and
JSPS Core-to-Core Program 18005.
\end{ack}


\end{document}